\documentclass{article}
\usepackage{amssymb}

\input{tcilatex}
\begin{document}

\title{\textbf{Inverse Laplace Transform for Bi-Complex Variables}}
\author{A. BANERJEE$^{1}$, S. K. DATTA$^{2}$ and MD. A.\ HOQUE$^{3}$ \\
%EndAName
$^{1}$Department of Mathematics, Krishnath College, Berhampore, \\
\ Murshidabad 742101, India, E-mail: abhijit.banerjee.81@gmail.com\\
$^{2}$Department of Mathematics, University of Kalyani, Kalyani, Nadia,\\
\ PIN-741235,India, E-mail: sanjib\_kr\_datta@yahoo.co.in \ \ \\
$^{3}$Gobargara High Madrasah (H. S.), Hariharpara, Murshidabad \\
\ 742166, India,E-mail: mhoque3@gmail.com}
\date{}
\maketitle

\begin{abstract}
In this paper we examine the existence of bicomplexified inverse Laplace
transform \ as an extension of it's complexified inverse version within the
region of convergence of bicomplex Laplace transform. In this course we use
the idempotent representation of bicomplex-valued functions as projections
on the auxiliary complex spaces of the components of bicomplex numbers along
two orthogonal,idempotent hyperbolic directions.

\end{abstract}

\section{Introduction}

The theory of bicomplex numbers is a matter of active research for quite a
long time science the seminal work of Segre[1] in search of special
algebra.The algebra of bicomplex numbers are widely used in the literature
as it becomes a valiable commutative alternative [2] to the non-commutative
skew field of quaternions introduced by Hamilton [3] (both are four-
dimensional and generalization of complex numbers).

A bicomplex number is defined as $\ $%
\[
\ \xi =a_{0}+i_{1}a_{1}+i_{2}a_{2}+i_{1}i_{2}a_{3}, 
\]%
where \ $a_{0,}$ $a_{1},a_{2},a_{3}$ are real numbers, $%
i_{1}^{2}=i_{2}^{2}=-1$ and 
\[
\ i_{1}i_{2}=i_{2}i_{1},(i_{1}i_{2})^{2}=1. 
\]

The set of bicomplex numbers,complex numbers and real numbers are denoted by 
$C_{2},C_{1},$and $C_{0}$ respectively. $C_{2}$ becomes a Real Commutative
Algebra with identity 
\[
\ \ 1=1+i_{1}\cdot 0+i_{2}\cdot 0+i_{1}i_{2}\cdot 0 
\]%
$\ \ \ \ \ $ with standard binary composition.

There are \ two non trivial elements $e_{1}=\frac{1+i_{1}i_{2}}{2}$ and $%
e_{2}=\frac{1-i_{1}i_{2}}{2}$ \ in $C_{2}$\ with the properties $\
e_{1}^{2}=e_{1},e_{2}^{2}=e_{2},e_{1}\cdot e_{2}=e_{2}\cdot e_{1}=0$ and $\
e_{1}+e_{2}=1$ which means that $\ e_{1}$ and $e_{2}$ are idempotents (some
times called also orthogonal idempotents). By the help of the idempotent
elements $\ e_{1}$ and $e_{2}$ any bicomplex number%
\[
\xi
=a_{0}+i_{1}a_{1}+i_{2}a_{2}+i_{1}i_{2}a_{3}=(a_{0}+i_{1}a_{1})+i_{2}(a_{2}+i_{1}a_{3})=z_{1}+i_{2}z_{2} 
\]%
\ where \ $a_{0,}$ $a_{1},a_{2},a_{3}\epsilon R,$ 
\[
z_{1}(=a_{0}+i_{1}a_{1}),z_{2}(=a_{2}+i_{1}a_{3})\epsilon C_{1} 
\]%
can be expressed as $\ $%
\[
\xi =z_{1}+i_{2}z_{2}=\xi _{1}e_{1}+\xi _{2}e_{2} 
\]%
where $\xi _{1}(=z_{1}-i_{1}z_{2})$ and $\xi _{2}(=z_{1}+i_{1}z_{2})\epsilon
C_{1}.$

This representation of a bicomplex number is known as the Idempotent
Representation of $\xi $. $\xi _{1}$ and $\xi _{2}$ are called the
Idempotent Components of the bicomplex number $\xi =z_{1}+i_{2}z_{2}$,
resulting a pair of mutually complementary projections 
\[
\ P_{1}:(z_{1}+i_{2}z_{2})\epsilon C_{2}\longmapsto
(z_{1}-i_{1}z_{2})\epsilon C_{1} 
\]
and 
\[
P_{2}:(z_{1}+i_{2}z_{2})\epsilon C_{2}\longmapsto (z_{1}+i_{1}z_{2})\epsilon
C_{1}. 
\]%
The spaces $A_{1}=\{P_{1}(\xi ):\xi \epsilon C_{2}\}$ and $A_{2}=\{P_{2}(\xi
):\xi \epsilon C_{2}\}$ are called the auxiliary complex spaces of bicomplex
numbers.

An element $\xi =z_{1}+i_{2}z_{2}$ is singular if and only if $%
|z_{1}^{2}+z_{2}^{2}|=0$.The set of singular elements is denoted as $O_{2}$
and defined by $O_{2}=\{\xi \epsilon C_{2}:\xi $ is the collection of all
-complex multiples of$\ e_{1}$ and $e_{2}$ $\}$

The norm the $||\cdot ||:C_{2}\longmapsto C_{0}^{+}$(set of all non negetive
real numbers) of a bicomplex number is defined as

\[
||\xi ||=\sqrt{\{|z_{1}|^{2}+|z_{2}|^{2}\}}=\sqrt{%
a_{0}^{2}+a_{1}^{2}+a_{2}^{2}+a_{3}^{2}} 
\]

\section{Laplace transform}

Let $f(t)$ be a real valued function of exponential order k. The \ coplex
version of Laplace Transform \ [5] of $f(t)$ for $t\geq 0$ can be defined as%
\[
L\{f(t)\}=F_{1}(\xi _{1})=\dint\limits_{0}^{\infty }f(t)e^{-\xi _{1}t}dt 
\]%
. Here $F_{1}(\xi _{1}:\xi _{1}\epsilon C_{1})$ exists and absolutely
convergent for Re$(\xi _{1})$ \TEXTsymbol{>} $k.$Similarly 
\[
F_{2}(\xi _{2})=\dint\limits_{0}^{\infty }f(t)e^{-\xi _{2}t}dt 
\]%
\ converges absolutely for Re$(\xi _{2})$ \TEXTsymbol{>} $k$\ \ .\ \ Then
the bicomplex Laplace Transform [4] of $f(t)$ for $t\geq 0$ can be defined as%
\[
L\{f(t)\}=F(\xi )=\dint\limits_{0}^{\infty }f(t)e^{-\xi t}dt 
\]%
. Here $F(\xi )$ exists and convergent in the region 
\[
D=\ \{\xi \epsilon C_{2}:\xi =\xi _{1}e_{1}+\xi _{2}e_{2}:Re(\xi
_{1})>k,Re(\xi _{2})>k\} 
\]
in idempotent representation.

\section{Inverse Laplace Transform for Bicomplex variables}

If $f(t)$ real valued function of exponential order k, defined on $t\geq 0$
,its Laplace \ transform $F_{1}(\xi _{1})$ in bicomplex variable $\xi
_{1}=x_{1}+i_{1}y_{1}\epsilon C_{_{1}}$ is simply%
\begin{eqnarray*}
F_{1}(\xi _{1}) &=&\dint\limits_{0}^{\infty }f(t)e^{-\xi
_{1}t}dt=\dint\limits_{0}^{\infty
}f(t)e^{-(x_{1}+i_{1}y_{1})t}dt=\dint\limits_{0}^{\infty
}e^{-x_{1}t}f(t)e^{-i_{1}y_{1}t}dt \\
&=&\dint\limits_{0}^{\infty
}\{e^{-x_{1}t}f(t)\}e^{-i_{1}y_{1}t}dt=\dint\limits_{-\infty }^{\infty
}g(t)e^{-i_{1}y_{1}t}dt=\psi (x_{1,}y_{1})
\end{eqnarray*}

which is Fourier transform of $g(t)$ where 
\[
g(t)=f(t)e^{-x_{1}t},t\geq 0;and=0,t<0 
\]
in usual complex exponential form.

$F_{1}(\xi _{1})$ converges for Re$(\xi _{1})$ \TEXTsymbol{>} $k$\ \ and

\[
|\ F_{1}(\xi _{1})|\ <\infty \ \Rightarrow |\dint\limits_{0}^{\infty
}f(t)e^{-\xi _{1}t}dt|=\dint\limits_{-\infty }^{\infty
}|g(t)e^{-i_{1}y_{1}t}|dt=\dint\limits_{-\infty }^{\infty }|g(t)|dt<\infty 
\]

The later condition shows that $g(t)$\ is absotulely integrable .Then by
Laplace inverse transform in complex exponential form

\[
g(t)=\frac{1}{2\pi i_{1}}\dint\limits_{-\infty }^{\infty
}e^{i_{1}y_{1}t}\psi (x_{1,}y_{1})dy_{1}\Rightarrow f(t)=\frac{1}{2\pi i_{1}}%
\dint\limits_{-\infty }^{\infty }e^{x_{1}t}e^{i_{1}y_{1}t}\psi
(x_{1,}y_{1})dy_{1}. 
\]%
\ \ \ \ \ \ \ \ \ \ \ \ \ \ \ \ \ \ \ \ \ \ \ \ \ \ \ \ \ \ \ \ \ \ \ \ \ \
\ \ \ \ \ \ \ \ \ \ \ \ \ \ \ \ \ \ \ \ \ \ \ \ \ \ \ \ \ \ \ \ \ \ \ \ \ \
\ \ \ \ \ \ \ \ \ \ \ \ \ \ \ \ \ \ \ \ \ \ \ \ \ \ \ \ \ \ \ \ \ \ \ \ \ 

\bigskip

Now if\ we integrate along a vertical line then x$_{1}$ is a constant and so
for complex variable $\xi _{1}=x_{1}+i_{1}y_{1}\epsilon C_{_{1}}($that
implies $d\xi _{1}=$ $dy_{1})$ \ the above inversion formula can be\ \ \ \ \
\ \ \ \ \ \ \ \ \ \ \ \ \ \ \ \ \ \ \ \ \ \ \ \ \ \ \ \ \ \ \ \ \ \ \ \ \ \
\ \ \ \ \ \ \ \ \ \ \ \ \ \ \ \ \ \ \ \ \ \ \ \ \ \ \ \ \ \ \ \ \ \ \ \ \ \
\ \ \ \ \ \ \ \ \ \ \ \ \ \ \ \ \ \ \ \ \ \ \ \ \ \ \ \ \ \ \ \ \ \ \ \ \ \
\ \ \ \ \ \ 

extended to complex Laplace inverse transform

\begin{eqnarray*}
f(t) &=&\frac{1}{2\pi i_{1}}\dint\limits_{x_{1}-i_{1}\infty
}^{x_{1}+i_{1}\infty }e^{(x_{1}+i_{1}y_{1})t}\psi (x_{1,}y_{1})dy_{1}=\frac{1%
}{2\pi i_{1}}\dint\limits_{x_{1}-i_{1}\infty }^{x_{1}+i_{1}\infty }e^{\xi
_{1}t}\psi (x_{1,}y_{1})d\xi _{1} \\
&=&\frac{1}{2\pi i_{1}}\lim\limits_{y_{1}\rightarrow \infty
}\dint\limits_{x_{1}-i_{1}y_{1}}^{x_{1}+i_{1}y_{1}}e^{\xi _{1}t}F(\xi
_{1})d\xi _{1}...........(1)
\end{eqnarray*}

$\ \ \ \ $

Here the integration is to be performed along a vertical line in the complex 
$\ \xi _{1}$-plane employing contour integration method.

We assume that $F_{1}(\xi _{1})$ \ is holomorphic in $x_{1}<k$ except for
having a finite number of poles $\xi _{1}^{k}$ ,$k=1,2,3,................n$
therein. Taking \ $R\rightarrow \infty $ we can guarantee

that all these poles lie inside the contour $\Gamma _{R}$ .Since $e^{\xi
_{1}t}$ never vanishes so the poles of $e^{\xi _{1}t}F(\xi _{1})$ and $%
F_{1}(\xi _{1})$ \ are same.Then by Cauchy residue theorem

\[
\lim\limits_{R\rightarrow \infty }\dint\limits_{\Gamma _{R}}e^{\xi
_{1}t}F(\xi _{1})d\xi _{1}=2\pi i_{1}\sum \func{Re}s\{e^{\xi _{1}t}F(\xi
_{1}):\xi _{1}=\xi _{1}^{k}\}. 
\]

Now since for $\xi $ on $C_{R}$ and $|F(\xi )$ $|<\frac{M}{|\xi |^{p}}$ [6]
some $p>0$ and\ all $R>R_{0},$%
\[
\ \lim\limits_{R\rightarrow \infty }\dint\limits_{C_{R}}e^{\xi _{1}t}F(\xi
_{1})d\xi _{1}=0\text{ }for\text{ }t>0\ \ \ 
\]%
.

so%
\[
\dint\limits_{\Gamma _{R}}e^{\xi _{1}t}F(\xi _{1})d\xi
_{1}=\dint\limits_{C_{R}}e^{\xi _{1}t}F(\xi _{1})d\xi
_{1}+\dint\limits_{x_{1}-i_{1}R}^{x_{1}+i_{1}R}e^{\xi _{1}t}F(\xi _{1})d\xi
_{1}=2\pi i_{1}\sum \func{Re}s\{e^{\xi _{1}t}F(\xi _{1}):\xi _{1}=\xi
_{1}^{k}\} 
\]

$\ \ \ \ \ \ \ \ \ \ \ \ \ \ \ \ \ \ \ \ \ \ \ \ \ \ \ .$

then for $R\rightarrow \infty $ we obtain

\[
\dint\limits_{x_{1}-i_{1}\infty }^{x_{1}+i_{1}\infty }e^{\xi _{1}t}F(\xi
_{1})d\xi _{1}=2\pi i_{1}\sum \func{Re}s\{e^{\xi _{1}t}F(\xi _{1}):\xi
_{1}=\xi _{1}^{k}\},t>0. 
\]

We first attend the right half plane D$_{1}=$Re$(\xi _{1})$ \TEXTsymbol{>} $%
k $ and

\[
\lim_{Re(\xi _{1})\longrightarrow \infty }F_{1}(\xi _{1})=0. 
\]
The inverse Laplace transform of $F_{1}(\xi _{1})$ will then a real valued
function

\[
f(t)=\frac{1}{2\pi i_{1}}\dint\limits_{x_{1}-i_{1}\infty
}^{x_{1}+i_{1}\infty }e^{\xi _{1}t}F_{1}(\xi _{1})d\xi _{1}\ \
..............\ \ \ \ \ (2) 
\]%
where$\ \xi _{1}=x_{1}+i_{1}y_{1}\epsilon C_{_{1}}.$

In the right half plane \ D$_{2}=$ Re$(\xi _{2})$ \TEXTsymbol{>} $k$ and

\[
\lim_{Re(\xi _{2})\longrightarrow \infty }F_{2}(\xi _{2})=0 
\]%
the inverse Laplace transform of $F_{2}(\xi _{2})$ will be

\[
f(t)=\frac{1}{2\pi i_{1}}\dint\limits_{x_{2}-i_{1}\infty
}^{x_{2}+i_{1}\infty }e^{\xi _{2}t}F_{2}(\xi _{2})d\xi _{2},\xi
_{2}=x_{2}+i_{1}y_{2}\epsilon C_{_{1}}..............(3) 
\]

Moreover in each case $f(t)$ \ is of exponential order k.

Then

\bigskip 
\begin{eqnarray*}
f(t) &=&f(t)e_{1}+f(t)e_{2}=\frac{1}{2\pi i_{1}}\dint\limits_{D_{1}}e^{\xi
_{1}t}F_{1}(\xi _{1})d\xi _{1}e_{1}+\frac{1}{2\pi i_{1}}\dint%
\limits_{D_{2}}e^{\xi _{2}t}F_{2}(\xi _{2})d\xi _{2}e_{2} \\
&=&\frac{1}{2\pi i_{1}}\dint\limits_{D=D_{1}\cup D_{2}}e^{\xi t}F(\xi )d\xi
..........(4)
\end{eqnarray*}%
where we use the fact that any real number $c$ can be written as 
\[
c=c+i_{1}\cdot 0+i_{2}\cdot 0+i_{1}i_{2}\cdot 0=c_{1}e_{1}+c_{2}e_{2}. 
\]%
The bicomplex version of inverse Laplace transform thus can be defined as
(4). Evidently, here also%
\[
\lim_{Re(\xi _{1,2})\longrightarrow \infty }F(\xi )=0 
\]%
and $f(t)$ is of exponential order $k$ . Reversing this proces one can at
once obtain f(t) \ from the integration defined in (4). It guarantees the
existance of inverse Laplace transform.

\subsection{Definition}

If $F(\xi )$ exists and is convergent in a region $D=D_{1}\cup D_{2}$ which
are the right half planes $D_{1,2}=R(\xi _{1,2})>k$ together with%
\[
\lim_{Re(\xi _{1,2})\longrightarrow \infty }F(\xi )=0 
\]
then the inverse Laplace transform of $F(\xi )$ can be defined as

\[
L^{-1}\{F(\xi )\}=\frac{1}{2\pi i_{1}}\dint\limits_{D=D_{1}\cup D_{2}}e^{\xi
t}F(\xi )d\xi =f(t) 
\]%
The integral in each plane $D_{1}$and $D_{2}$ are taken along any straight
line $R(\xi _{1,2})>k$ . As a result our object function $f(t)$ will be of
exponential order $k$ ,in the principal value sense.

\subsection{Examples}

\begin{itemize}
\item If we take $F(\xi )d\xi =\frac{1}{\xi }$, then it's inverse Laplace
transform is given by

\[
f(t)=\frac{1}{2\pi i_{1}}\dint\limits_{D=D_{1}\cup D_{2}}e^{\xi t}F(\xi
)d\xi =\frac{1}{2\pi i_{1}}\dint\limits_{D_{1}}e^{\xi _{1}t}F_{1}(\xi
_{1})d\xi _{1}e_{1}+\frac{1}{2\pi i_{1}}\dint\limits_{D_{2}}e^{\xi
_{2}t}F_{2}(\xi _{2})d\xi _{2}e_{2}............(4) 
\]

$\ \ \ \ \ \ $

Now 
\[
\frac{1}{2\pi i_{1}}\dint\limits_{D_{1}}e^{\xi _{1}t}F_{1}(\xi _{1})d\xi
_{1}=\frac{1}{2\pi i_{1}}\dint\limits_{x_{1}-i_{1\infty }}^{x_{1}+i_{1\infty
}}e^{\xi _{1}t}\frac{1}{\xi _{1}}d\xi _{1}=2\pi i_{1}\cdot 1=2\pi i_{1} 
\]

as $\xi _{1}=0$ is the only singular point therein, \ so%
\[
residue=\lim_{\xi _{1}\longrightarrow 0}(\xi -0)e^{\xi _{1}t}\frac{1}{\xi
_{1}}=1. 
\]
\end{itemize}

In a similar way,%
\[
\frac{1}{2\pi i_{1}}\dint\limits_{D_{2}}e^{\xi _{2}t}F_{2}(\xi _{2})d\xi
_{2}=2\pi i_{1} 
\]
\ and those leads (4) to%
\[
f(t)=e_{1}+e_{2}=1. 
\]

\begin{itemize}
\item In our procedure one may easily check a partial list....to name a
few....

\item $L^{-1}\{\frac{\omega }{\xi ^{2}+\omega ^{2}}\}=\sin \omega t,$

\item $L^{-1}\{\frac{\xi }{\xi ^{2}+\omega ^{2}}\}=\cos \omega t,$

\item $L^{-1}\{\frac{\xi +a}{(\xi +a)^{2}+\omega ^{2}}\}=e^{-at}\cos \omega
t,$

\item $L^{-1}\{\frac{\omega }{(\xi +a)^{2}+\omega ^{2}}\}=e^{-at}\sin \omega
t.$
\end{itemize}

References:

[1]C.Segre Math.Ann.:40,1892,pp:413

[2] N.Spampinato Ann.Math.Pura .Appl.:14,1936,pp:305

[3] W.R .Hamilton Lectures on quaternion Dublin:Hodges and Smith:1853

[4] A.Kumar,P.Kumar IJET : 3,2011,pp225.

[5] Y.V.Sidorov,M.V.Fedoryuk,M.I.Shabunin Mir Publishers,Moscow: 1985.

[6] Joel L.Schiff \ Springer.

\end{document}